\documentclass[a4paper,oneside,11pt]{article}

\usepackage{SarahEnglHead}
\usepackage{undertilde}
\usepackage{Commandes}
\usepackage[all]{xy}
\usepackage{epic}
\usepackage{graphicx}
\usepackage{float}
\restylefloat{figure}
\title{Existence of Natural and Projectively Equivariant Quantizations}
\author{Sarah Hansoul\\
 Department of Mathematics\\ 
University of Li\`{e}ge\\
Grande Traverse, 12, B-4000 Li\`{e}ge\\
Belgium\\
email: \texttt{s.hansoul@ulg.ac.be}
}
\begin{document}
\maketitle

{Keywords: Quantizations, Projective Structures, Differential Operators, Smooth Manifolds}

\medskip

{MSC [2000]: 58A05, 53B05, 46L65, 58J70}

\begin{abstract}
We study the existence of natural and projectively equivariant quantizations for differential operators acting between order $1$ vector bundles over a smooth manifold $\M$. To that aim, we make use of the Thomas-Whitehead approach of projective structures and construct a Casimir operator depending on a projective Cartan connection. We attach a scalar parameter to every space of differential operators, and prove the existence of a quantization except when this parameter belongs to a discrete set of resonant values. 
\end{abstract}

\section{Introduction}

The term quantization carries several meanings in mathematics. It seems to be an appropriate starting point for this introduction to define what will be meant by quantization throughout this paper.
Given two natural vector bundles $E_1$ and $E_2$ over a smooth manifold $\M$, one can consider the space $\mathcal{D}(E_1,E_2)$ of differential operators acting between smooth sections of $E_1$ and $E_2$. This space is filtered by the order of differentiation:

\[
\mathcal{D}(E_1,E_2)=\bigcup_{k\in\N}\mathcal{D}^k(E_1,E_2).
\]
The space $\mathcal{S}(E_1,E_2)$ of {principal symbols} associated to $\mathcal{D}(E_1,E_2)$ is the graded space associated to this filtration:

\[
\mathcal{S}(E_1,E_2)=\bigoplus_{k\in\N}\mathcal{S}^k(E_1,E_2),
\]
where $\mathcal{S}^k(E_1,E_2)=\mathcal{D}^k(E_1,E_2)/\mathcal{D}^{k-1}(E_1,E_2)$.

It is immediate to see that $\mathcal{D}(E_1,E_2)$ and $\mathcal{S}(E_1,E_2)$ are isomorphic as vector spaces. Moreover, each of these spaces is 
a representation of the group of local diffeomorphisms. One could thus wonder if they remain isomorphic once endowed with this representation structure.  
If we call \textit{quantization} a linear bijection
\[
Q:\mathcal{S}(E_1,E_2)\to\mathcal{D}(E_1,E_2),
\]
satisfying a normalization condition, then the latter question is equivalent to asking for the existence of a \textit{natural} quantization, namely commuting with the action of local diffeomorphisms. An infinitesimal version of naturality is the commutation of the quantization with Lie derivatives in the direction of any element of the algebra $\Vect(\M)$ of smooth vector fields over $\M$.


It is not difficult to see that such quantizations do not exist. Though, some intermediary results have been obtained in the Euclidian case $\M=\R^m$ (see \cite{duval1, duval2, duval3, lecomte3, lecomte4, lecomte5}).

In \cite{lecomte4}, P. Lecomte shows the existence and the unicity of an $\slmp$-invariant quantization for the space $\mathcal{D}_{\lambda\mu}$ of differential operators acting between $\lambda$ and $\mu$ densities, when the difference $\delta=\mu-\lambda$ does not belong to a discrete set of so-called resonnant values. An $\slmp$-invariant quantization is a quantization $Q$ which satisfies
\[
\Lie{X}\circ Q=Q\circ\Lie{X}\pourtt X\in\slmp,
\]
where $\slmp$ is the {projective} subalgebra of $\Vect(\R^m)$. It is the embedding of $\slmu$ generated by the vector fields
\[
\partial_i, x^j\partial_i,\mbox{ and } x^jx^i\partial_i.
\]
Thus, though $\mathcal{D}_{\lambda\mu}$ is not isomorphic to its space of symbols as a representation of $\Vect(\R^m)$, it is isomorphic to it as a representation of $\slmp$. This result is in a sense optimal since $\slmp$ is a {maximal} subalgebra of the Lie algebra $\Vect_*(\R^m)$ of polynomial vector fields over $\R^m$.

Similar results have been obtained on $\R^m$ for the spaces $\mathcal{D}_{\lambda\mu}$: in \cite{duval2}, C. Duval, P. Lecomte and V. Ovsienko show the existence and unicity of a $so_{p+1,q+1}$-invariant quantization (with $p+q=m$), where $so_{p+1,q+1}$ is the conformal algebra of $\Vect(\R^m)$, also maximal in $\Vect_*(\R^m)$. More generally, F. Boniver and P. Mathonet have determined all the maximal subalgebras $\g$ of $\Vect_*(\R^m)$ and settled an existence criterion for $\g$-invariant quantizations (\cite{boniver2}).

Given the abundance of results obtained over $\R^m$, the question naturally arises whether or not they can be extended to an arbitrary manifold $\M$. In this case, one can not ask for $\slmp$-invariance of a quantization since $\slmp$ is not well defined anymore: there is no canonical embedding of $\slmu$ in $\Vect(\M)$ in general.  P. Lecomte proposed the below generalization of $\slmp$-invariant quantizations in \cite{lecomte2}, using the notion of projective structure. A projective structure (or projective class) over a manifold $\M$ is an equivalence class of linear torsion free connections on $\M$ defining the same geodesics over $\M$ up to reparametrization.

Let $\Co{\M}$ denote the space of linear torsion free connections over $\M$. A \textit{natural and projectively equivariant quantization} corresponding to a space $\mathcal{D}(E_1,E_2)$ of differential operators is a map 
\[
Q:\Co{\M}\times\mathcal{S}(E_1, E_2)\to\mathcal{D}(E_1,E_2),
\]
such that
\begin{itemize}
\item[$\bullet$] $Q$ is \textit{natural}, \emph{i.e.} commutes with the action of local diffeomorphisms;
\item[$\bullet$] $Q$ is \textit{projectively equivariant}: $Q_\nabla=Q_{\nabla'}$ if $\nabla$ and $\nabla'$ are projectively equivalent;
\item[$\bullet$] $Q_\nabla:\mathcal{S}(E_1,E_2)\to\mathcal{D}(E_1,E_2)$ is a quantization for all $\nabla\in\Co{\M}$.
\end{itemize}

Defined this way, natural and projectively equivariant quantizations are a generalization of the notion of $\slmp$-invariant quantizations over $\R^m$ in the following sense: if $\nabla^0$ denotes the canonical flat connection over $\R^m$ and $Q$ is such as above, then $Q_{\nabla^0}$ is a $\slmp$-invariant quantization.

In this paper, we settle an existence criterion for natural and projectively equivariant quantizations when $E_1$ and $E_2$ are vector bundles associated to the fiber bundle $\PM$ of linear frames of $\M$.

Our method relies on M. Bordemann's proof of the existence of natural and projectively equivariant quantizations for differential operators acting between densities (\cite{Bord}). The method consists in translating the problem from the original manifold $\M$ to a manifold $\Mt$, which is a principal bundle over $\M$ with one additionnal dimension. The choice of the bundle $\Mt$ is justified by the existence of a natural map associating a connection $\nablat$ to any projective class $\cl{\nabla}$ of linear connections over $\M$. This construction corresponds to the approach undertaken by T.Y. Thomas and J.H.C. Whitehead in the '$20$s of the geometry of paths, which is the study of projective structures (see \cite{thomas, veblenthom, whitehead}).

The quantization problem is easily solved on $\Mt$ by means of the so-called {standard quantization} $\tilde{\tau}$, a natural map providing a quantization for any linear connection over $\Mt$. It is in addition projectively equivariant when considered as depending on a linear connection $\nabla$ over $\M$ via $\nablat$.

Using this method, the difficulty lies in the translation of the problem from $\M$ to $\Mt$. The first step consists in transforming a symbol over $\M$ into a symbol over $\Mt$ by a natural and projectively equivariant map $\Rel$. After having applied the standard quantization to the lifted symbol, one has to bring down the differential operator obtained over $\Mt$ on a differential operator over $\M$, still via a natural and projectively equivariant map $\Des$. The situation is summarized by the following diagram
\begin{equation}\label{diagramme}
\xymatrix@=.2em{
\tilde{\mathcal{S}}({E_1},{E_2})\ar[rrrr]^{\taut}&&&&
\tilde{\mathcal{D}}({E_1},{E_2})\ar[dddd]^{\Des_\nabla}\\
&&&&\\
&&&&\\
&&&&\\
\mathcal{S}(E_1,E_2)\ar[uuuu]^{\Reln}&&&&
\mathcal{D}(E_1,E_2)\\
}
\end{equation}
The searched quantization is given by $Q_\nabla=\Des_\nabla\circ\taut\circ\Reln$ for every $\nabla$ in $\Co{\M}$.

A first attempt to adapt M. Bordemann's results  to differential operators acting between $p$-forms and smooth functions showed that the lift of symbols is a very critical step in this construction (\cite{pformes}). In the case of differential operators acting between densities, the lifted symbols were solutions of a partial derivatives equation. Long and tedious computations were necessary to find the equivalent of this equation for differential operators acting on forms. Moreover, there was no clue for a systematic construction of such an equation for arbitrary differential operators.

The first problem we encounter in the search of a general lift of symbols is the determination of the space $\tilde{\mathcal{S}}(\E_1,\E_2)$ of symbols over $\Mt$ in which the lift was valued. We there use the assumption made on the vector bundles we consider: since they are associated 
to $\PM$, symbols are assimilated to functions from $\PM$ to $\V$, where $\V$ is a representation of $\GLm$. In section \ref{Vtilde}, we show how it is possible to associate a representation $\Vt$ of $\GLmu$ to any representation $\V$ of $\GLm$. The lift of symbols which are functions from $\PM$ to $\V$ is defined as valued in functions from $\PMt$ to $\Vt$.

To tackle the lift problem, we use the theory of Cartan bundles and Cartan connections. Its interest in our problem hinges on the fact that a Cartan connection $\w_{\nabla}$ on a Cartan bundle over $\M$ is naturally associated to any projective structure $\cl{\nabla}$ over $\M$. Moreover, we define a Casimir operator over any Cartan bundle, whose definition and properties are similar to the properties of the classical Casimir operator acting on a representation of a semi-simple Lie algebra but which depends on a Cartan connection. This operator plays a key role in our construction of the lift of symbols. To perform it, we do not make use anymore of a partial derivatives equation but rather define the lift of a symbol $P$ to be the unique eigenvector of the Casimir operator depending on $\w_{\nabla}$ which projects on $P$.

A scalar parameter $\delta$ is attached to every space of differential operators in a way described in section \ref{Vtilde}. When it belongs to a discrete set of resonant values, it is not possible to lift the symbols by means of the Casimir operator. We explicitely compute the set of these values for every space of differential operators, and in the case of $\mathcal{D}_{\lambda\mu}$ recover  the set of resonant values for $\delta=\mu-\lambda$ obtained over $\R^m$ in \cite{lecomte4}. On the other hand, the descent of differential operators is quite straightforward to perform for any value of $\delta$.

Altogether, we obtain the existence of a natural and projectively equivariant quantization for differential operators acting between sections of order $1$ vector bundles whenever $\delta$ is not resonnant. Note that this result is not optimal: P. Mathonet and F. Radoux have shown in the case of differential operators acting between densities that there are resonant values for which the quantization still exists (\cite{pierrefabian}). These values are the ones for which the $\slmp$-quantization exists over $\R^m$ but is not unique.

There are still open questions related to the existence of natural quantizations. First, $\slmp$-invariant quantizations are unique for non resonant values. One could wonder if it is still the case for natural and projectively equivariant quantizations; F. Radoux  showed that unicity is not preserved (personal communication). Thus, one could classify these quantizations for example by a cohomology approach, initiated in \cite{affcohom}.

Second, other invariances than the projective one could be generalized over an arbitrary manifold. In particular one could study the existence of natural and \emph{conformally} equivariant quantizations, namely natural quantizations depending on a metric in a conformally equivariant way. We have the hope that an existence result could be obtained via a method similar to the one presented in this paper, but using the bundle described by C. Fefferman and C. Graham in \cite{Feffgra} rather than $\Mt$.

\section{Basic notions and notations}\label{Basic}
Throughout this paper, we will denote by $\M$ a smooth, Hausdorff, connex, countable manifold. Moreover, we will assume that the dimension $m$ of $\M$ is strictly greater than $1$.

\medskip

\subsection{Differential operators and associated symbols}

Let $E_1(\M)$, $E_2(\M)$ (or simply $E_1$, $\E_2$ when no confusion is possible) be two finite rank vector bundles over $\M$, and $\mathcal{D}(E_1,E_2)$ the space of differential operators acting between the sets $\Gamma^\infty(E_1)$ and $\Gamma^\infty(E_2)$ of smooth sections of $E_1$ and $E_2$ respectively. This space is filtered by the order of differentiation: 
\[
\mathcal{D}(E_1,E_2)=\bigcup_{k\in\N}\mathcal{D}^k(E_1,E_2),
\]
where $\mathcal{D}^k(E_1,E_2)$ is the space of differential operators of order at most $k$. Locally, such operators are of the form 
\[
\mathcal{D}=\sum_{\va\alpha \leq
k}C_\alpha\partial_1^{\alpha_1}\cdots\partial_m^{\alpha_m},
\]
where the coefficients $C_{\alpha}$ belong to $\mathrm{Hom}(E_1,E_2)$.

\medskip

\begin{defi}The \textit{principal symbol} map, denoted by $\sigma$, is the map defined (in local coordinates) on $\mathcal{D}(\E_1,\E_2)$ by
\[
\sigma: 
\sum_{\va\alpha \leq
k}C_\alpha\partial_1^{\alpha_1}\cdots\partial_m^{\alpha_m}\mapsto\sum_{\va\alpha =
k}C_\alpha\xi_1^{\alpha_1}\cdots\xi_m^{\alpha_m}.
\]
\end{defi}
We will denote by $\mathcal{S}(E_1,E_2)$ the space of \textit{principal symbols} (or simply \textit{symbols}) of the elements of $\mathcal{D}(E_1,E_2)$. This space is nothing else but the graded space associated with the filtration of the latter:
\[
\mathcal{S}(E_1,E_2)=\bigoplus_{k\in\N}\mathcal{S}^k(E_1,E_2),
\]
where each $\mathcal{S}^k(E_1,E_2)$ is isomorphic to $\Gamma^{\infty}(\mathrm{Hom}(E_1,E_2)\otimes\sym^kT\M)$.

\medskip 
\begin{defi}
A \textit{quantization} for $\mathcal{D}(\E_1,\E_2)$ is a linear bijection
\[
\mathcal{Q}: \mathcal{S}(E_1,E_2)\to\mathcal{D}(E_1,E_2)
\] 
satisfying $\sigma\circ \mathcal{Q}=\id$.

\end{defi}
Such maps are easy to find. The \textit{standard quantization}, or \textit{Lichnerowicz quantization}, is a map
\[
\tau:\Co{\M}\times\mathcal{S}(E_1,E_2)\to\mathcal{D}(E_1,E_2),
\]
where $\Co{\M}$ denotes the space of torsion free linear connections over $\M$ (see \cite{lichne}).
This map provides a quantization for every connection $\nabla$ and is natural: for every local diffeomorphism $\phi$ over $\M$, one has
\[
\phipb(\tau_\nabla(P))=\tau_{(\phipb\nabla)}(\phipb P)\pourtt\nabla\in\Co{\M}, P\in \mathcal{S}(\E_1,\E_2).
\]

 The existence of quantizations implies that every space of differential operators is isomorphic as vector space to its space of symbols. Though, they don't remain isomorphic once endowed with their canonical structure of representation of local diffeomorphisms. In other words, among the quantizations between $\mathcal{S}(E_1,E_2)$ and $\mathcal{D}(E_1,E_2)$, none is {natural}, \emph{i.e.} satisfies for every local diffeomorphism $\phi$ over $\M$ and every symbol $P$ $\phipb(\mathcal{Q}(P))=\mathcal{Q}(\phipb P)$.

\subsection{Order $1$ bundles}
In this paper, we won't consider arbitrary vector bundles over $\M$.  Recall that if $\G$ is a Lie group, $P$ is a $\G$-principal bundle over $\M$ and $F$ is a manifold on which $\G$ acts on the left, then the	\textit{associated bundle} to $P$  and $F$, denoted by $\ass{P}{\G}{F}$
is the quotient of $P\times F$ by the equivalence relation $\sim$ defined by
\[
(p,f)\sim(p',f') \mbox{ if there exists } g\in\G \mbox{ such that }(p',f')=(p\cdot g,g^{-1}\cdot f).
\]
It is a classical result that sections of these bundles are in canonical correspondence with $\G$-equivariant functions $h:P\to F$, namely functions satisfying
\[
h(x\cdot g)=g^{-1}\cdot h(x)\pourtt g\in\G. 
\] 
\begin{defi}
Let $\PM$ be the $\GLm$-principal bundle of linear frames of $\M$:
\[
\PM=\{\mbox{frames } v_x\mbox{ of }T_x\M: x\in\M\}.
\] 
A vector bundle $E$ over $\M$ is said to be \textit{of order $1$} if it is associated to $\PM$:
\[
E=\ass{\PM}{\GLm}{\V}
\]
for some (left) representation $(\V,\rho)$ of $\GLm$.
\end{defi}

In this paper, we will always assume $\E_1$ and $\E_2$ to be of  order $1$, covering tensor and densities bundles as a particular case. We will systematically assimilate sections of these bundles to their associated equivariant functions over $\PM$. Given a representation $(\V,\ro)$ of $\GLm$, the space of these functions will be denoted by 
\[
\finv{\PM}{\V}{\ro}.
\] 
As a consequence of this choice, the symbols of differential operators will also be sections of an order $1$ bundle. Indeed, if $\E_1$ and $\E_2$ are associated to $\PM$ and respectively $\V_1$ and $\V_2$, then the homomorphisms from $\E_1$ to $\E_2$ correspond to equivariant functions
 from $\PM$ to $\V_1^*\otimes\V_2$. Thus any symbol of order $k$ is assimilated to an element of 
\[
\finv{\PM}{\V_1^*\otimes\V_2\otimes\sym^k\R^m}{\ro}.
\]

\medskip
\begin{rem}
An important feature of the vector bundles belonging to this family is their naturality: any local diffeomorphism $\phi:\M\to N$ can be lifted into a local diffeomorphism $\phipb:\ass{\PM}{\GLm}{\V}\to\ass{P^1N}{\GLm}{\V}$. Therefore, it makes sense to talk about the naturality of the maps we will consider in the sequel.
\end{rem}

\subsection{The question}
\begin{defi}
Two connections $\nabla$ and $\nabla'$ in $\Co{\M}$ are said to be \textit{projectively equivalent} if there exists $\alpha\in\Omega^1(\M)$ such that
\[
\nabla'_X Y=\nabla_X Y+\frac{1}{m+1}(\alpha(X)Y+\alpha(Y)X)\pourtt
X,Y\in\Vect(\M),
\]
where $\Vect(\M)$ is the Lie algebra of smooth vector fields over $\M$. Alternatively, H. Weyl proved that two connections are projectively equivalent if and only if they define the same geodesics on $\M$ up to reparametrization (\cite{weyl}). The equivalence classes of this relation are called \textit{projective structures}, or \textit{projective classes}, and the projective class of a given connection $\nabla$ will be denoted $\cl{\nabla}$.
\end{defi}
 \vspace{1em}
 \begin{defi}A \textit{natural and projectively equivariant quantization} for $\mathcal{D}(\E_1,\E_2)$ is an operator associating to every manifold $\M$ a map  
\[
Q^\M:\Co{\M}\times\mathcal{S}(E_1(\M),
E_2(\M))\to\mathcal{D}(E_1(\M),E_2(\M)),
\]
such that\\
\begin{itemize}
\item[$\bullet$] $Q^\M$ is \textit{natural}: for all local diffeomorphism $\phi:\M\to N$ one has
\[
Q^\M_{(\phi^*\nabla)}(\phi^* P)=\phi^*(Q^N_\nabla(P)),
\]
for all $\nabla\in\Co{N}$ and $P\in \mathcal{S}(E_1(N),E_2(N))$;\\
\item[$\bullet$] $Q^M$ is \textit{projectively equivariant}
(or \textit{projectively invariant}): 
$Q^M_\nabla=Q^M_{\nabla'}$ if $\nabla$ and $\nabla'$ are projectively
equivalents;\\
\item[$\bullet$] $Q^M_\nabla:\mathcal{S}(E_1(\M),E_2(\M))\to\mathcal{D}(E_1(\M),E_2(\M))$ is a quantization for all $\nabla\in\Co{\M}$.
\end{itemize}
\end{defi}
\medskip
In the sequel to avoid uselessly heavy notations we will simply denote by $Q$ the map $Q^\M$.

\medskip

The aim of this paper is to give an existence criterion for natural and projectively equivariant quantizations when $E_1$ and $E_2$ are order $1$ vector bundles.

\section{Two approaches of projective structures}

In the 20's, several authors beared interest in the so-called projective geometry of path, namely the study of projective structures. While one approach, developed among others by V. Oblen, T.Y. Thomas and J.H.C. Whitehead, is formulated in terms of symmetric linear connections, E. Cartan independently established at almost the same time the theory of Cartan connections to study projective structures.  These theories are of course linked, but depending on the context, it can be more convenient to work with one of the two formalisms. As both of them have been useful to solve our problem, we briefly present each of them in this section.

\subsection{Thomas-Whitehead's approach}
At the time the original papers of T.Y. Thomas (\cite{thomas, veblenthom}) and J.H.C. Whitehead (\cite{whitehead}) were written, the theory of principal bundles and Ehresmann connections was not yet established. C.W. Roberts formulated their work in terms of this more modern setting in \cite{Roberts}. 
The idea behind this approach consists in associating a unique linear connection on a $m+1$-dimensional principal bundle over $\M$ to a projective structure on a $m$-dimensional manifold $\M$.

\medskip
We will denote by $\Mt\stackrel{\pi}{\to}\M$ the bundle $\ass{\PM}{\GLm}{\R_0^+}$, where the action of $\GLm$ over $\R_0^+$ is given by 
\[
A\cdot \lambda=\va{\det A}^{-1}\lambda.
\]
It is the fiber bundle of positive $1$-densities over $\M$. Remark that it is also a $\R_0^+$-principal bundle for the action 
\[
\nu:\Mt\times\R_0^+\to\Mt;\; (\cl{v,r},s)\mapsto\cl{v, r s}.
\]
Therefore, we can consider on $\Mt$ the \textit{Euler} vector field $\Eul$ which is the fondamental vector field associated to $1$: $\Eul_y=\dto{(y\cdot e^t)}$ for all $y\in\Mt$. 

\medskip

\begin{defi}
A \textit{Thomas-Whitehead projective connection}, or \textit{T-W connection} is a torsion-free linear connection $\hat{\nabla}$ over $\Mt$ satisfying
\begin{itemize}
\item[$\bullet$]$\hat\nabla\Eul=\frac{1}{m+1}\id$;
\item[$\bullet$]$\displaystyle{\nu_{s*}(\hat\nabla_X
Y)=\hat\nabla_{\nu_{s*}(X)}\nu_{s*}(Y)\pourtt X,Y\in\Vect(\Mt)}$.
\end{itemize}
\end{defi}
A particular family of T-W connections consists in the \emph{normal} T-W connections, whose curvature satisfies a normalization condition (see \cite{Roberts} for a complete definition).

Any T-W connection canonically induces a projective structure on $\M$ in a way described in \cite{Roberts}. Thomas and Whitehead's approach of projective structures lies on the following fundamental theorem, proved in \cite{Roberts}. 
\medskip
\begin{thm}
Given a projective structure $[\nabla]$ on a manifold $\M$, there exists a unique normal T-W connection $\nablat$ on $\Mt$ inducing $[\nabla]$. Moreover, the map $\nabla\mapsto\nablat$ is natural.
\end{thm}
\subsection{Cartan's approach}
While Thomas and Whitehead associate to a projective structure a unique normal T-W connection on the bundle $\Mt$, Cartan associates to a projective structure a normal projective Cartan connection on another principal bundle over $\M$ called {Cartan bundle}. We here only outline the theory of Cartan connections, a more detailed presentation can be found in \cite{Kob} or \cite{Sharpe} for example.
\medskip

\begin{defi}
Let $\M$ be a $m$-dimensional manifold, $G$ be a Lie group with Lie algebra $\Galg$, $G_0$ be a closed subgroup of $G$ such that $\dim G/G_0=m$ and $P$ be a $G_0$-principal bundle over $\M$. A \textit{Cartan connection} $\w$ over $P$ is a $\Galg$-valued $1$-form over $P$ satisfying
\begin{itemize}
\item[$\bullet$]$\w(\fond{h})=h \mbox{ for all } h$ in the Lie algebra $\Goalg$ of
$G_0$;
\item[$\bullet$]$\Right{g}^*\w=\Ad (g^{-1})\circ\w \mbox{ for all } g\in\Go$;
\item[$\bullet$]$\w_y$ is a bijection between $T_y\M$ and $\Galg$
 for all $ y\in P$,
\end{itemize}
where $\Right{}$ denotes the (right) action of $\Go$ over $P$.
\end{defi}
\vspace{1em}
\begin{defi}
The \textit{curvature} of a Cartan connection is the $\Galg$-valued $2$-form $\Omega$ over $P$ defined by
\[
\Omega=\mathrm{d}\w+\nij{\w,\w},
\]
where $\mathrm{d}$ is the De Rham differential and 
\[
\nij{\w,\w}(X,Y)=\nij{\w(X),\w(Y)}\pourtt X,Y\in TP.
\]
\end{defi}

As it is the case with T-W connections, there is a notion of normality for Cartan connections. \emph{Normal} Cartan connections are defined by means of a normalization condition that their curvature satisfies.

\medskip

\begin{defi}
A Cartan connection is said to be \textit{projective} if
\[
\G=\PGLm=\GLmu/\R_0,
\]
and
\begin{eqnarray*}
\Go=\Hmu&=&\{
\matgr{A}{0}{\alpha}{a}:A\in\GLm,a\in\R_0,\alpha\in\R^{m*}\}/\R_0\\
&=&\{\matHmu{A}{0}{\alpha}{1}:A\in\GLm,\alpha\in\R^{m*}\}.
\end{eqnarray*}
\end{defi}
The Lie algebras of $\G$ and $\Go$ are then respectively given by
\[
\Galg=\glmu/\R=\slmu,
\]
and
\[
\Goalg=\hmu=\{\matgr{A}{0}{\alpha}{-\tr A}:A\in\GLm,\alpha\in\R^{m*}\}.
\]
Note that the Lie algebra $\Galg$ is $3$-graded of type
\[
\Galg=\gmun\oplus\go\oplus\gun,
\]
where $\gmun=\{\matgr{0}{u}{0}{0}:u\in\R^{m}\}$ is isomorphic to $\R^m$, $ \go=\{\matgr{A}{0}{0}{-\tr A}:A\in\glm\}$
 is isomorphic to $\glm$ and
$\gun=\{\matgr{0}{0}{\alpha}{0}:\alpha\in\R^{m*}\}$ is isomorphic to $\R^{m*}$.

\subsection{A particular Cartan bundle}

A projective Cartan bundle is simply a bundle admitting a projective Cartan connection. The latter condition does not fix an explicit form for such a bundle but rather prescribes how it behaves under a change of coordinates (see \cite{Sharpe}). An explicit model of a Cartan bundle is for example given in \cite{Kob} but we here describe another Cartan bundle, proposed by D.J. Saunders and M. Crampin in \cite{Crampin}, which is closer to the original idea of Cartan of "attaching a projective space of the same dimension than the manifold to each of its points" (\cite{Cartan}).

\medskip
\begin{defi} For every $X\in T\Mt$, let $\cl{X}$ denote the class of $X$ for the equivalence relation induced by the action $\nu_*$ of $\R_0^+$ on $T\Mt$: 
\[
\cl{X}=\cl{Y} \mbox{ if there exists }s\in\R_0^+\mbox{ such that } X=\nu_{*s}Y.
\]
The bundle $\TM$ over $\M$ is the space of all these equivalence classes:
\[
\TM=\{\cl{X}: X\in T\Mt\}.
\]
\end{defi}

It is a vector bundle whose typical fiber is $m+1$-dimensional. Thus, by taking the quotient of each fiber by the multiplicative action of $\R_0$ we  obtain over each point of $\M$ a projective space of dimension $m$. Moreover, $\TM$ admits a global non vanishing section $e$, given by
\[
e_x= \cl{\Eul_y}, 
\]
for any $y$ in the fiber $\pi^{-1}(x)$ of $x$.

\medskip

Consider the frame bundle associated to  $\TM$. It is a $\GLmu$-principal bundle on which $\R_0$ acts canonically by multiplication.
Taking the quotient $\PrM$ of the latter by this action, we obtain a $\PGLm$-principal bundle over $\M$. 
The Cartan bundle we will use in our proof of the existence of quantizations is a subbundle of the latter. 

\medskip

\begin{defi}
The bundle $\CM$ is the $\Hmu$-principal subbundle of $\PrM$ containing the points whose last element is a (non zero) multiple of $e$: 
\[
\CM=\{(\cl{X_1}_x,\ldots,\cl{X_{m}}_x,\alpha e_x)\mbox{ basis of
}\mathcal{T}_x\M: x\in\M,\alpha\in\R_0\}/\R_0.
\]
This bundle can also be viewed as
\[
\{(\cl{X_1}_x,\ldots,\cl{X_{m}}_x,e_x)\mbox{ basis of }\mathcal{T}_x\M:
x\in\M\}.
\]

\end{defi}
\begin{nota}\label{w}${}$
In order to avoid lengthy notations, the points of $\CM$ will be denoted by $(\cl{w},e)$, where 
$w=(X_1,\ldots,X_m)$ and $\cl{w}=(\cl{X_1},\ldots,\cl{X_m})$. 
\end{nota}
The following result is proved in \cite{Crampin}.

\medskip
\begin{thm}
To every projective class $\cl{\nabla}$ on $\M$ is associated a unique normal Cartan connection $\w_\nabla$ over $\CM$. Moreover, the correspondence $\nabla\mapsto\w_\nabla$ is natural.
\end{thm}
\medskip
\begin{rem}
The two described approaches of projective structures are linked to each other by the theory of \emph{tractor calculus} (see \cite{capgov1, capgov2}). This general theory extends beyond the purpose of this paper but we have implicitely used it in the above developements: the bundle $\TM$ actually belongs to the family of \emph{tractor bundles}.
\end{rem}

\section{From a representation of $\GLm$ to a representation of $\GLmu$}\label{Vtilde}

Following the general proof strategy explained in the introduction and summarized by the diagram (\ref{diagramme}), we will have to lift elements of $\mathcal{S}(\E_1,\E_2)$ to a corresponding space $\tilde{\mathcal{S}}(\E_1,\E_2)$ over $\Mt$. But the intuitive notion of \textit{corresponding} space of symbols $\tilde{\mathcal{S}}(E_1,E_2)$ over $\Mt$ to a given space of symbols $\mathcal{S}(E_1,E_2)$ over $\M$ is not well defined. As an example, while it seems obvious to lift elements of $\mathcal{S}(C^\infty(\M),C^\infty(\M))$ to $\mathcal{S}(C^\infty(\Mt),C^\infty(\Mt))$, one could wonder whether the elements of $\mathcal{S}(\wedge^m(T^*M),C^\infty(\M))$ should be lifted to  $\mathcal{S}(\wedge^m(T^*\Mt),C^\infty(\Mt))$ or to $\mathcal{S}(\wedge^{m+1}(T^*\Mt),C^\infty(\Mt))$.

To answer this question and obtain a formal definition of the space of symbols over $\Mt$ in which the original symbols will be lifted, we are going to associate to every representation $\V$ of $\GLm$ a representation $\Vt$ of $\GLmu$. This is sufficient to tackle the above problem since we are only considering vector bundles associated to $\PM$. This implies that each element of 
$\mathcal{S}^k(\E_1,\E_2)$ is identified to a $\GLm$-equivariant function 
\[
\PM\to\V_1^*\otimes\V_2\otimes\sym^k\R^m.
\] 
We will then define its lift as valued in the space of $\GLmu$-equivariant functions
\[
\PMt\to(\V_1^*\otimes\V_2\otimes\sym^k\R^m)^\sim.
\]

This section is dedicated to the definition of $\Vt$ given $\V$, and to the study of some properties of the representations of the type $\Vt$. This is performed using the theory of representations of $\GLm$, that we here outline following \cite{fulhar}, \cite{goodwall} and \cite{Murnaghan}. We will only consider finite dimensional and continuous representations of $\GLm$, \textit{i.e.} representations $(\V,\rho)$ such that $\rho(A)$ is continuous in the components $(A)_{ij}$ of $A$ for every $A\in\GLm$.

\subsection{Young diagrams and irreducible representations}\label{young}

\begin{defi} A \textit{Young diagram}, or
\textit{Ferrer diagram}
 of \textit{size} $d\in\N$ and \textit{depth} $p\in\N$  is an element
  $\D=(d_1,\cdots,d_{p})$ of $\N^p$ satisfying\\
\[
d_1\geq\ldots\geq d_p,\quad\mbox{ and }\quad\sum_i d_i=d.
\]
\end{defi}
A Young diagram $\D=(d_1,\ldots,d_p)$ is commonly represented by an array of $p$ lines, the $i^{th}$ line containing $d_i$ boxes. The following picture represents the Young Diagram $(3,2,2)$.

\vspace{1em}

\begin{center}
\input{PictArt1.epic}
\end{center}
 
Young diagrams are of particular interest in representation theory since they allow a handy classification of irreducible representations of $\GLm$ (see \cite{Murnaghan}). 
\medskip

On the one hand, every $(\D,n,\delta)$ where $\D$ is a Young diagram of depth at most $m-1$, $n\in\Z$ and $\delta\in\C$, canonically gives rise to an irreducible representation $S_{(\D,n,\delta)}$ of $\GLm$. It is included in $\otimes^d\R^m$ if $\D=(d_1,\cdots,d_{m-1})$ is of size $d$, and given by
\[
S_{(\D,n,\delta)}=\{\va{\det A}^{\delta}(\det A)^n\roc(A)v_\D:A\in\GLm\},
\]
where $\roc$ denotes the canonical action of $\GLm$ over $\otimes^d\R^m$, and
\[
v_\D=e_1^{d_1-d_2}\otimes(e_1\wedge e_2)^{d_2-d_3}\otimes\cdots\otimes(e_1\wedge\cdots\wedge e_{m-1})^{d_{m-1}}.
\]

On the other hand, each irreducible representation of $\GLm$ is characterized up to isomorphism by a Young diagram of depth at most $m-1$, an integer and a complex number, \emph{i.e.} is isomorphic to some $S_{(\D,n,\delta)}$. 
\medskip

\begin{nota}
In the sequel, unless otherwise stated, $(\V,\ro)$ will denote an irreducible representation characterized by $(\D,n,0)$. A representation characterized by $(\D,n,\delta)$ will be denoted $(\V,\rodel)$.
\end{nota}

We will of course deal with non irreducible representations of $\GLm$ but we will always assume  them to be \textit{finite direct sums of irreducible representations of $\GLm$}. Moreover, the above classification allows us to consider without loss of generality that irreducible representations of $\GLm$ are always included in $\otimes^d\R^m$ for some $d\in\N$.

\subsection{Adding one dimension}

\begin{defi}
Let $(\V,\rodel)$ be an irreducible representation of $\GLm$ characterized by $(\D,n,\delta)$. Then the irreducible representation $(\Vt,\rot)$ of $\GLmu$ induced by $\V$ is the one characterized by $(D,n,0)$.
\end{defi}
If $\V$ is not irreducible but admits the decomposition $\V=\V_1\oplus\cdots\oplus\V_p$, the above definition is extended by 
\[
\Vt=\Vt_1\oplus\cdots\oplus\Vt_p.
\]

It will be of further relevance to study the behavior of the correspondence $\V\to\Vt$ relatively to the correspondence $\V\to\V^*$ between a representation and its dual. Once more, Young diagrams are a convenient tool to study this behaviour. If $\V$ is characterized by $(\D,n,\delta)$, then $\V^*$ is characterized by $(\D^*,-n-d_1,-\delta)$, where $d_1$ is the length of the first line of $\D$ and $\D^*$ is obtained  by taking the complement of $\D$ in a diagram of $m$ lines of $d_1$ boxes.

This implies that the representation $\Vete$ is not isomorphic to $\Vt$: it is characterized by the diagram $\ete{\D}$  obtained by adding above $\D$ a line of the same length as the first line of $\D$, as illustrated in the figure below. This is due to the fact that $\D^*$ is complementary to $\D$ according to a diagram of depth $m$, while the $\ete{\D}$ is complementary to $\D^*$ according to a diagram of depth $m+1$.
\begin{center}
\includegraphics[scale=0.4]{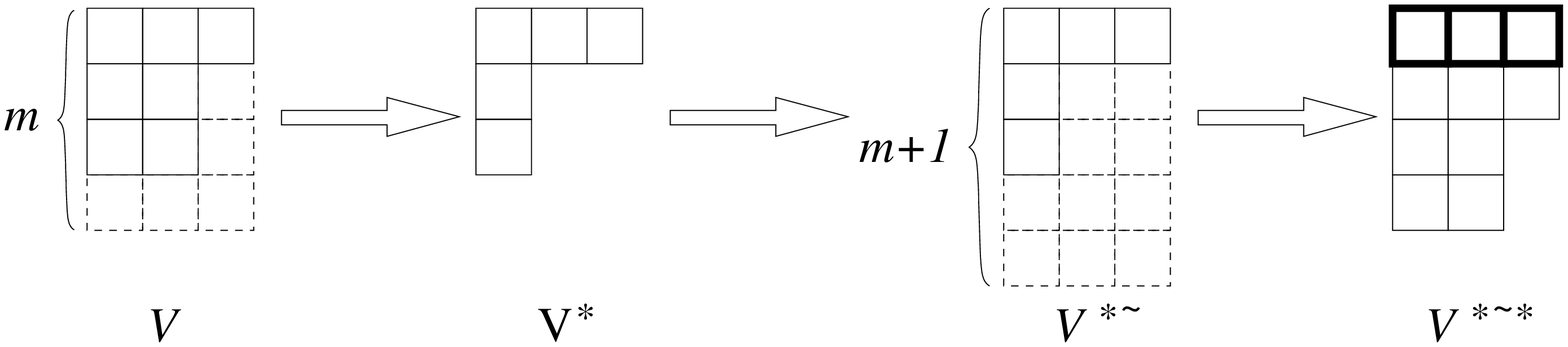}
\end{center}
\subsection{Branching law}
Whenever $(\V,\ro)$ is an irreducible representation of $\GLm$ $(\Vt,\rot)$ is, by definition, an irreducible representation of $\GLmu$. But it is also a representation of $\GLm$ for 
\[
\rep(A)=\rot\matgr{A}{0}{0}{1}.
\]
Endowed with the latter structure, $\Vt$ is not irreducible in general but admits a well known decomposition described by the following proposition, known as \textit{branching law} (see \cite{goodwall}).

\vspace{1em}
\begin{prop}\label{branch}
Let $(W,\rep)$ be a irreducible representation of $\GLmu$, characterized by $(\D=(d_1,\ldots,d_m),n,\delta)$. Then its irreducible components under $\GLm$ are exactly the ones characterized by $(E=(e_1,\ldots,e_{m-1}),n,\delta)$, where
\[
d_1\geq e_1\geq\cdots\geq d_{m-1}\geq e_{m-1}\geq d_m.
\]
Each of these components is of multiplicity $1$ in the decomposition.
\end{prop}
Graphically, the branching law states that the Young diagrams that describe the irreducible components of a representation are obtained by removing boxes to the original diagram, respecting the rule that a line in a new diagram can not be shorter than the lower line in the original diagram, as shown in the example below.

\begin{figure}[H]
\begin{center}
\includegraphics[scale=0.4]{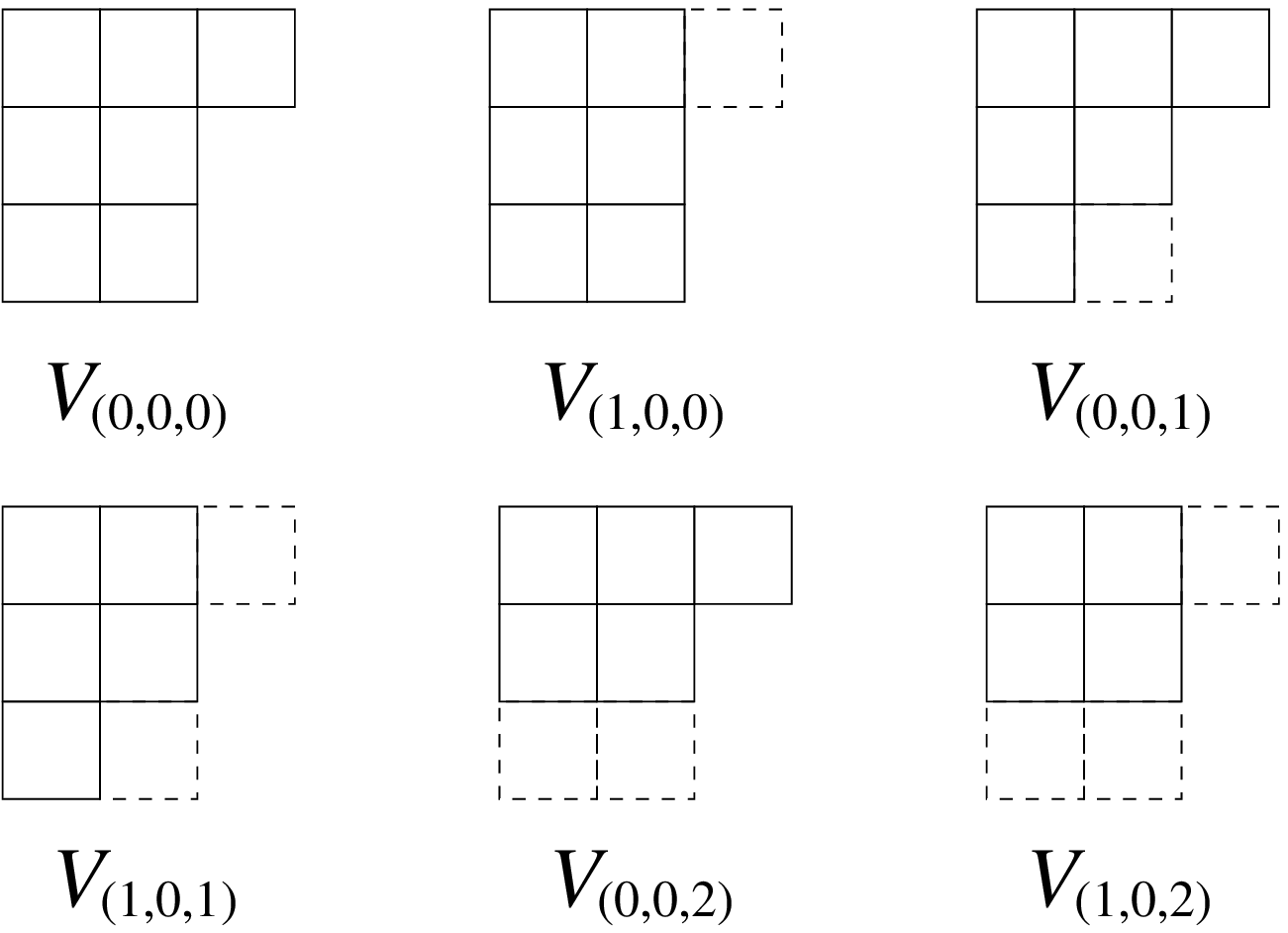}
\caption{Young diagrams appearing in the decomposition of  $\Vt$ when $\V$ corresponds to the diagram $(3,2,2)$.}
\end{center}
\end{figure}

\begin{nota}\label{vvqq}
We will adopt the notation used in the above figure: every irreducible component of $\Vt$ will be labelled by some $q\in\N^{m-1}$ and noted $\V_q$. The $i^{\scriptscriptstyle{th}}$ component $q_i$ of $q$ represents the number of boxes removed from the original diagram at line $i$ to obtain the diagram describing $\V_q$. The set of all possible labels of the irreducible components under $\GLm$ of $\Vt$ will be denoted $Q$ instead of $Q_{\Vt}$ when no confusion is possible.
\end{nota}

The branching law admits the following interpretation: as a representation of $\GLmu$,  $\R^{m+1}$ is characterized by the diagram consisting in one single box. It decomposes into $\R^m\oplus\R$ as shown below. 

\begin{center}
\includegraphics[scale=0.4]{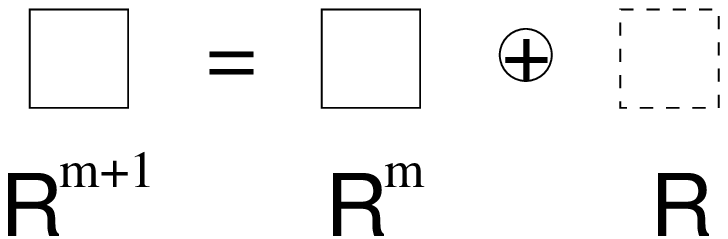}
\end{center}

\noindent This corresponds to the decomposition
\[
\vect{X}{x}=\vect{X}{0}+\vect{0}{x}
\]
of any point in $\R^{m+1}$. The decomposition described by the branching law is simply the extension of the latter decomposition to $\otimes^d\R^{m+1}$. In a given subrepresentation of $\Vt$, every removed box from the original diagram corresponds to a term of the form $\vect{0}{x}$ while the remaining boxes correspond to terms of the form $\vect{X}{0}$. 

\medskip
\begin{rem}\label{subrep}
As illustrated in Figure $1$, a representation $\V$ of $\GLm$ is always isomorphic to a subrepresentation of $\Vt$: it corresponds to the diagram where no boxes have been removed, namely to the irreductible component $\V_0$ of $\Vt$. With the above interpretation, $V_0$ is the component whose elements do not contain terms of the form $\vect{0}{x}$.

A similar phenomenon occurs when considering $\Vete$.  Indeed, by removing as many boxes as possible according to the branching law to the diagram characterizing $\ete{\V}$, one always recover the diagram describing $\V$ as illustrated in the figure below.
\begin{figure}[H]
\begin{center}
\includegraphics[scale=0.4]{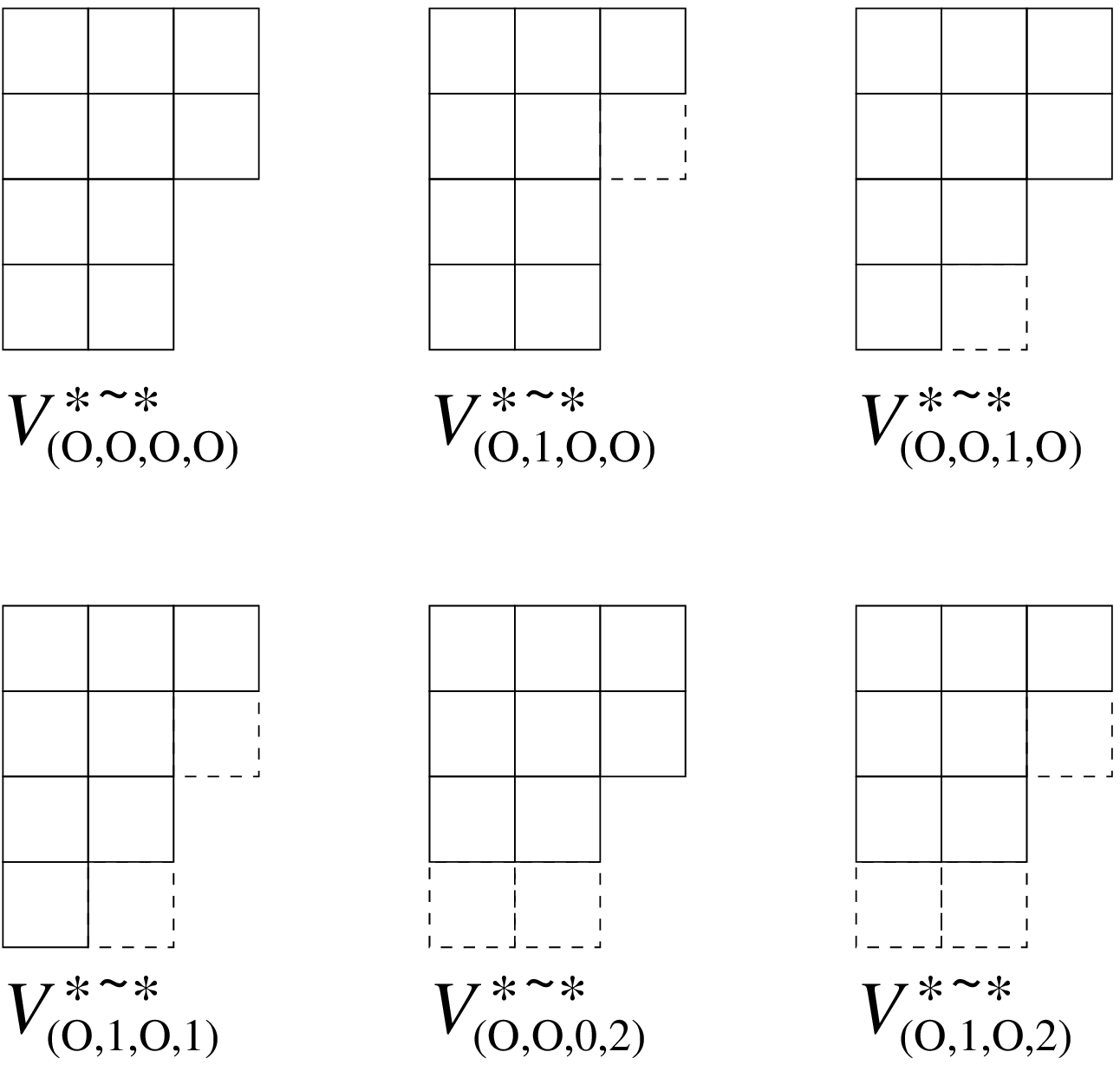}
\end{center}
\caption{Young diagrams appearing in the decomposition of $\Vete$ when $\V$ corresponds to the diagram $(3,2,2)$.}
\end{figure}
The original representation $\V$ is nested in $\Vete$ but this time it is isomorphic to the irreducible component containing the maximum number of terms of the form $\vect{0}{x}$.
\end{rem}
\begin{nota}
If $\E$ denotes an order $1$ vector bundle corresponding to a representation $\V$ of $\GLm$, then we will denote $\tilde{\E}$ and $\ete{E}$ the vector bundles associated to $\PMt$ and to $\Vt$ and $\Vete$ respectively.
\end{nota}

\subsection{Resonant values}
In the discussion of the existence of the {$\slmp$-invariant quantization} for differential operators acting between $\lambda$ and $\mu$ densities over $\R^m$ (\cite{lecomte4}), several values of the shift $\mu-\lambda$, called \textit{resonant values} naturally arise, for which the quantization is not unique or does not exist. Since the existence of natural and projectively equivariant quantizations as considered in this paper implies the existence of $\slmp$-invariant quantization over $\R^m$, these values also appear over an arbitrary manifold. The phenomenon of resonant values is not specific to the case of differential operators acting between densities. We describe thereafter how every irreducible representation $(\V,\ro)$ of $\GLm$ carries a finite set of resonant values.

\vspace{1em}

Recall that the \textit{Casimir operator} on a representation $(W,r)$ of a semi-simple  $p$-dimensional Lie algebra $\g$ acts on $W$ by 
\[
\cas{}=\sum_i r(u_i)\circ r(u_i^+),
\] 
where $(u_1,\ldots,u_p)$ and $(u_1^+,\ldots,u_p^+)$ are two Killing-dual basis of $\g$.
\vspace{1em}

The semi-simple Lie algebra we will consider here is the embedding $\slmp$ of $\slmu$ in the Lie algebra $\Vect(\R^m)$ of smooth vector fields over $\R^m$. It is the maximal subalgebra of the Lie algebra of polynomial vector fields over $\R^m$ spanned by
\[
\partial_i,x^i\partial_j\mbox{ and }x^ix^k\partial_k.
\]

Every representation $\V$ of $\GLm$ induces a structure of $\slmp$-representation on the sections of $\ass{P^1\R^m}{\GLm}{\V}$, via the Lie derivative
\[
\Lie{X}s=X\cdot s+\rho_*(DX)\circ s\pourtt X\in\slmp,
\]
where $DX$ denotes the differential of the vector field $X$. Thus, a Casimir operator can be defined on this space of sections. Such Casimir operators have been computed in \cite{boniver3}, where we find the following result.  
\medskip
\begin{prop}
Let  $\V$ be an irreducible representation of $\GLm$ characterized by $(\D,n,\delta)$. Then the Casimir operator of $\slmp$
 on $\Gamma^{\infty}(\ass{\PM}{\GLm}{\V})$ is equal to a multiple of the identity. If $\D$ is equal to $(d_1,\ldots,d_{m-1})$ and is of size $d$, then this multiple is equal to 

\begin{multline}\label{alpha}
\al=\frac{(m(n-\delta)+d)(m(n+1-\delta)+d)}{2m}\\
+\frac{1}{2m(m+1)}\sum_{i,j=1}^m
d_id_j(m\delta_{ij}-1)+2d_i(m-j)(m\delta_{ij}-1).
\end{multline}
\end{prop}

We will call \textit{eigenvalue associated to $(\V,\rodel)$} the number $\alpha$ above defined.
\medskip
\begin{defi} 
Let $(\V,\rho_\epsilon)$ be an irreducible representation of $\GLm$. A complex number $\delta$ is \textit{resonant} for $\V$ if, given that $(\Vt,\rot_{\delta+\epsilon})$ decomposes into $\bigoplus_{q\in Q}\V_q$, there exists $q\in Q\backslash\{0\}$ such that $\al_q= \al_0$, where $\al_q$ denotes the eigenvalue associated to $\V_q$ and $\al_0$ the eigenvalue associated to $\V_0$.

If $\V$ is not irreducible, $\delta$ is said to be resonant for $\V$ if it is resonant for one of its irreducible components.
\end{defi}
Given the explicit form of the eigenvalues given in (\ref{alpha}), it is straightforward to compute the set of resonant values for irreducible representations of $\GLm$.
\medskip
\begin{prop}
The set of resonant values for a representation characterized by 
$(\D,n,0)$ is
\[
\bigl\{
\frac{\va{q}(2(m+1)(n+1)+2d-\va{q})+\sum_{i=1}^{m}2d_iq_i-q_i^2-2iq_i}{2\va{q}(m+1)}:q\in
Q\backslash \{0\} \bigr\},
\]
where $\va{q}=\sum_iq_i$ and $\D=(d_1,\ldots,d_{m-1})$ is of size $d$.
In particular, there exists only a finite number of resonant values for any representation $(\V,\ro)$ of $\GLm$, and $0$ is never resonant whenever $n\geq 0$.
\end{prop}

\section{Construction of the quantization}
As  explained in the introduction, the main obstacle to adapt M. Bordemann's proof of the existence of natural and projectively equivariant quantizations lies in the lift of symbols. With the tools described in the previous section at hand,  it is possible for every representation $\V$ of $\GLm$ to lift in a natural and projectively equivariant way elements of $\finv{\PM}{\V}{\rodel}$ into elements of $\finv{\PMt}{\Vt}{\rot}$. Since we are working with order $1$ fiber bundles, lift of symbols is a particular case of this construction. Applying standard quantization $\taut$ to the lifted symbols, there remains to bring down from $\Mt$ to $\M$ the differential operators obtained to find the searched quantization. The descent of differential operators appears surprisingly straightforward. 

\subsection{Casimir operator depending on a Cartan connection}
In the proof of the existence of the $\slmp$-invariant quantization over $\R^m$, a key role is held by the Casimir operator relative to $\slmp$ acting both on symbols and on differential operators. But this operator can not be defined on an arbitrary manifold, because $\slmp$ is not canonically nested anymore in $\Vect(\M)$.

On the other hand, every projective Cartan connection $\w$ on $\CM$ gives rise to an embedding of $\slmu$ in $\Vect(\M)$, via $\w^{-1}$. In a broader context, every $\Galg$-valued Cartan connection over a bundle $\fib{P}$ allows to embed $\Galg$ in $\Vect(\fib{P})$. We can therefore define the  \textit{Casimir operator} depending on a Cartan connection, so called since its definition is similar to the definition of the classical Casimir operator on representations of semi-simple Lie algebras.  
\medskip
\begin{defi}
Let $\Go$ be a closed subgroup of a Lie group $\G$, $\w$ be a Cartan connection over a $\Go$-principal bundle $\fib{P}\to\M$ such that
$\Galg$ is a $p$-dimensional semi-simple Lie algebra, and let $(u_1,\ldots,u_{p})$ and $(u^+_1,\ldots,u^+_{p})$ be two Killing-dual basis of $\Galg$. Then the \textit{Casimir operator} relative to $\w$ is defined by
\[
\cas{\w}=\sum_{i=1}^{p}\Lie{\w^{-1}(u_i)}\circ \Lie{\w^{-1}(u^+_i)}.
\]
\end{defi}
This definition is independent on the choice of the basis $(u_1,\ldots,u_{p})$ of $\Galg$.
\medskip

\begin{prop}\label{natcas}
The Casimir operator above defined is natural.
\end{prop}

\begin{dem}We have to prove that for any local diffeomorphism $\phi$ over $P$,
every Cartan connection $\w$ over $\fib{P}$ and every smooth function $f$ over $P$, 
one has
\[
\phipb(\casw f)=\cas{\phipb\w}\phipb f.
\]

To that aim, it is sufficient to remark that for all $k\in\Galg$ the vector fields $(\phipb\w)^{-1}(k)$ and $\w^{-1}(k)$ are $\phi$-linked. Therefore
\[
\phipb(\Liew{k}f)=\Lie{((\phipb\w)^{-1}k)}(\phipb f),
\]
and in turn
\begin{eqnarray*}
\phipb(\casw f)&=&\phipb(\sum_{i=1}^{p}\Lie{\w^{-1}(u_i)}\circ
\Lie{\w^{-1}(u^+_i)}f)\\
&=&\sum_{i=1}^{p}\Lie{(\phipb\w)^{-1}(u_i)}\circ \Lie{(\phipb\w)^{-1}(u^+_i)}(\phipb f)\\
&=&\cas{\phipb\w}(\phipb f). 
\end{eqnarray*}
\end{dem}
\medskip
\begin{rem}
The usual Casimir operator on a representation of $\Galg$ has the property to commute with the action of the elements of $\Galg$. It is not longer the case with the Casimir operator depending on a Cartan connection, though the $\Ad$-invariance of the latter implies the following weaker property: if $\w$ is a Cartan connection, then for every $h\in\Goalg$, one has
\[
\Liew{h}\circ\casw=\casw\circ\Liew{h}.
\]
\end{rem}

The Casimir operator defined above can be computed explicitely on two families of smooth functions over the projective Cartan bundle $\CM$. The first one is the following.

Consider a representation $(\V,\ro)$ of $\GLm$. It is endowed with a structure $(\V,t)$ of representation of $\Hmu$ defined by
\[
t(\matHmu{A}{0}{\alpha}{1})=\ro(A).
\]
We will first consider the set 
\[
\finv{\CM}{\V}{\sdel}
\]
of smooth functions from $\CM$ to $\V$ satisfying $f(u\cdot\matHmu{A}{0}{\alpha}{1})=\va{\det A}^{\delta}t(\matHmu{A}{0}{\alpha}{1}^{-1})f(u)$ for all $\matHmu{A}{0}{\alpha}{1}$ in $\Hmu$ and $u$ in $\CM$. Using the symmetries of these functions and properties of Cartan connections, a direct computation leads to the following result.
\medskip
\begin{prop}Let $f$ be an element of $\finv{\CM}{\V}{\sdel}$, and let $\al$ be the eigenvalue associated to $\V$. Then for any Cartan connection $\w$ on $\CM$,
\[
\casw f=\al f.
\]
\end{prop}

Similarly, if $(\Vt,\rot)$ is a representation of $\GLmu$, it is also endowed with a structure $(\Vt,\tilde{t})$ of representation of $\Hmu$:

\[
\tilde{t}(\matHmu{A}{0}{\alpha}{1})=\rot(\matgr{A}{0}{\alpha}{1}).
\]
As above, we denote
\[
\finv{\CM}{\Vt}{\tdel}
\]
the set of smooth functions from $\CM$ to $\Vt$ satisfying $g(u\cdot\matHmu{A}{0}{\alpha}{1})=\va{\det A}^{\delta}\tilde{t}(\matHmu{A}{0}{\alpha}{1}^{-1})g(u)$ for all $\matHmu{A}{0}{\alpha}{1}$ in $\Hmu$ and $u$ in $\CM$. As a prolongation to Notation \ref{vvqq}, we will denote by $g_q$ the projection of $g$ on the irreducible component $\V_q$ of $\Vt$ for every $q\in Q$.

It is also possible to explicitely compute the Casimir operator on these functions.
\medskip
\begin{prop}\label{castinv}
 Let $g\in\finv{\CM}{\Vt}{\tdel}$. Then
 \[
\casw g=\sum_{q\in Q}\al_q g_q-2\sum_{q\in Q}\sum_{i=1}^m
i(\epsilon^i)\circ\Liew{v_i} g_q,
 \]
where $\al_q$ is the eigenvalue associated to $\V_q$ for every $q\in Q$, $v_i$ is equal to $\frac{2}{m+1}\matgr{0}{e_i}{0}{0}$ and $i(\epsilon^i)$ is defined on $\otimes^d\R^{m+1}$ by
\[
i(\epsilon^i)(\vect{X_1}{x_1}\otimes\ldots\otimes\vect{X_d}{x_d})=\sum_{j=1}^d X^i_j\vect{X_1}{x_1}\otimes\cdots\otimes\vect{0}{x_j}\otimes\cdots\otimes\vect{X_d}{x_d}.
\]
\end{prop}
Restricted to any irreducible component $\V_q$, the Casimir operator on $\finv{\CM}{\Vt}{\tdel}$ is thus the sum of a multiple of the identity and of a term valued in every $\V_{q'}$ such that $\va{q'}=\va{q}+1$, namely whose elements have one more term of the form $\vect{0}{x}$ than the elements of $\V_q$.

\vspace{1em}

\subsection{Lift of symbols}
The symbols we are considering in this paper are always elements of $\finv{\PM}{\V}{\rodel}$ for some representation $(\V,\rodel)$ of $\GLm$. In order to obtain an easy push down of differential operators, we will ask their lift to be valued in $\finv{\PMt}{\Vt}{\rot}^{-\delta}$. The latter is the subspace of $\finv{\PMt}{\Vt}{\rot}$ whose elements satisfy
\[
g(\nu_{s*}v)=s^{-\delta}g(v),\;\forall s\in\R_0^+.
\]  
Such functions will be reffered as $(-\delta)$-equivariant functions.

The Casimir operator will reveal to be a convenient tool, but it is only defined on $\CM$; we therefore transform the original lift problem into an equivalent lift problem over $\CM$. Following the below diagram, the new lift problem consists in finding a natural map $\mathcal{R}$ depending on a Cartan connection.

\begin{equation}
\xymatrix@=.2em{\finv{\PMt}{\Vt}{\rot}^{-\delta}&&&&
\finv{\CM}{\Vt}{\tdel}\ar[llll]^{\tilde{u}}\\
&&&&\\
&&&&\\
&&&&\\
\finv{\PM}{\V}{\rodel}\ar[rrrr]^{u}&&&& \finv{\CM}{\V}{\sdel}\ar[uuuu]^\Rw}
\end{equation}

Altogether, the searched lift will be given by
\[
\Reln=\tilde{u}\circ\mathcal{R}_{\w_\nabla}\circ u.
\]

\begin{lem}\label{run}

There exists a natural bijection
\[
u:\finv{\PM}{\V}{\ro_{\delta}}\to\finv{\CM}{\V}{\sdel}.
\]
\end{lem}
\begin{dem}
The map $u$ can be defined using the explicit form of the Cartan bundle we are working with. Using Notations \ref{w}, if $(\cl{w},e)$ is a point of $\CM$ then $\pi_*w$ is a frame of $\M$. Thus, we can set
\[
u(f)(\cl{w},e)=f(\pi_* w)\pourtt f\in\finv{\PM}{\V}{\rodel}.
\]
The map $u(f)$ is well defined due to the $\rodel$-equivariance of $f$. The bijectivity and naturality of $u$ are immediate to check.
\end{dem}

\begin{lem}\label{rdeux}
There exists a natural bijection
\[
\tilde{u}:\finv{\CM}{\Vt}{\tdel}\to\finv{\PMt}{\Vt}{\rot}^{-\delta}.
\]
\end{lem}
\begin{dem}
This time, we have to use the fact that $\Mt$ is associated to $\PM$. It implies that for any $(\cl{w},e)$ in $\CM$ and any $y$ in $\Mt$, there exists a unique $r\in\R_0^+$ such that $y=\cl{\pi_*w,r}$. It allows us to define $\tilde{u}$ by
\[
\tilde{u}(g)(w,\Eul)_{\cl{\pi_*w,r}}=r^{-\delta}g(\cl{w},e)\pourtt g\in\finv{\CM}{\Vt}{\tdel}.
\]
The above only defines $\tilde{u}(g)$ on frames of $\Mt$ of the form $(w,\Eul)$, but the equivariance of $g$ allows extension of $\tilde{u}(g)$ on the whole $\PMt$ by $\rot$-equivariance. The term $r^{-\delta}$ makes $\tilde{u}$ valued in $(-\delta)$-equivariant functions. Again, the bijectivity and naturality of this map are easy to check.
\end{dem}
\medskip
These two lemmas ensure us that the lift problem will be solved if we can find a natural map depending on a Cartan connection from $\finv{\CM}{\V}{\sdel}$ to $\finv{\CM}{\Vt}{\tdel}$. Note that every $\tdel$-equivariant function $g$ valued in $\Vt$ canonically gives rise to a $\sdel$-equivariant function valued in $\V$: the latter is nothing else than its projection $g_0$ on the component $\V_0$ of $\Vt$. Though, this correspondence is not injective: there are several ways of prolongating a $\tdel$-equivariant function whose projection on $\V_0$ is known. A natural and projectively equivariant one is given by the following proposition whose proof, quite long, is given at the end of the paper.
\medskip 
\begin{prop}\label{rw}
Let $\V$ be a representation of $\GLm$ and $\delta\in\R$ be non resonant for $\V$. Given $f\in\finv{\CM}{\V}{\sdel}$, there
exists a unique function $\Rw(f)\in\finv{\CM}{\Vt}{\tdel}$ such that
\begin{itemize}
\item[$\bullet$] $\Rw(f)$ is an eigenvector of $\casw$;
\item[$\bullet$] $(\Rw(f))_0=f$.
\end{itemize}
Moreover, the map $\mathcal{R}$ is natural. 
\end{prop}

As announced, the lift existence is a direct corollary of the three above results. 
 \medskip
\begin{thm}\label{Reln} If $\V$ is a representation of $\GLm$ and $\delta\in\R$ is not resonant for $\V$, then there exists a natural and projectively equivariant map
\[
\Rel:\Co{\M}\times\finv{\PM}{\V}{\rodel}\to\finv{\PMt}{\Vt}{\rot}^{-\delta}.
\]
It is given by $\Reln=\tilde{u}\circ\mathcal{R}_{\w_\nabla}\circ u$.
\end{thm}

\subsection{Descent of differential operators}
When working with symbols of differential operators in $\mathcal{D}(\E_1,\E_2)$, one could expect their lift to be the symbol of a differential operator in $\mathcal{D}(\tilde{\E}_1,\tilde{\E}_2)$.
This is not the case however. We know that any symbol of order $k$ over $\M$, which is valued in $\V_1^*\otimes\V_2\otimes\sym^k\R^{m}$, is lifted into a function valued in $(\V_1^*\otimes\V_2\otimes\sym^k\R^{m})^\sim$. It is a consequence of Littlewood-Richardson's rule, whose statement extends beyond the purpose of this paper and can be found in \cite{Murnaghan} that $\widetilde{\V\otimes\V'}\subset\Vt\otimes\Vt'$
 for any representations $\V$ and $\Vt'$ of $\GLm$. In consequence, the lifted symbols are valued in $\V_1^{*\sim}\otimes\Vt_2\otimes\sym^k\R^{m+1},$
which means they are symbols of differential operators in $\mathcal{D}(\ete{\E_1},\tilde{\E}_2)$, acting from  $\finv{\PMt}{\Vete_1}{\rot_1}^{-\lambda}$ to $\finv{\PMt}{\Vt_2}{\rot_2}^{-\mu}$.

\medskip

This observation makes our task easier: due to the way representations $\V$ are embedded in $\Vete$, functions in $\finv{\PM}{\V}{\rodel}$ are actually in one-to-one correspondence with functions in $\finv{\PMt}{\Vt}{\rot}^{-\delta}$.

\medskip
\begin{lem}\label{eteu}
For every representation $(\V,\rodel)$ of $\GLm$, there exists a natural bijection
\[
 \ete{u}:\finv{\PM}{\V}{\ro_{\delta}}\to\finv{\PMt}{\Vete}{\rot}^{-\delta}.
\]
\end{lem}
\begin{dem} We consider the case where $(\V,\ro)$ is an irreducible representation of $\GLm$, $\ete{u}$ being defined component wise in the general case.

Remark \ref{subrep} tells us that $\V$ is always isomorphic to an irreducible component $\Vete_p$ of $\Vete$. Thus, any $\V$-valued function over $\PM$ can be viewed as a $\Vete$-valued map. Let us define $\ete{u}$ by
\[
\ete{u}(f)(w,\Eul)_{\cl{\pi_* w,r}}=r^{-\delta}f(\pi_*w)\pourtt f\in\finv{\PM}{\V}{\rodel}.
\]
This definition implies that every function $\ete{u}(f)$ has the same values on two frames $(w,\Eul)$ and $(w',\Eul)$ where $\pi_*w=\pi_*w'$, in other words such that $(w',\Eul)=(w,\Eul)\cdot\matgr{\Id}{0}{\alpha}{1}$ for some $\alpha\in\R^{m*}$. We have to check that the latter is compatible with the $\rot$-equivariance of $\ete{u}(f)$, which implies
\begin{eqnarray*}
\ete{u}(f)(w',\Eul)&=&\rot(\matgr{\Id}{0}{\alpha}{1})^{-1}\ete{u}(f)(w,\Eul).
\end{eqnarray*}
Due to the way $\GLmu$ acts on $\otimes^d\R^{m+1}$, the right hand side of this equality is equal to the sum of $\ete{u}(f)(w,\Eul)$ and of other tensors containing at least one more term of the form $\vect{0}{x}$ than the latter, which belongs by definition to $\Vete_p$. Since we've seen that $\Vete_p$ is the component of $\Vete$ whose elements contain the maximum number of such terms, these additional tensors have to vanish.

This point being checked, every $\ete{u}(f)$ is extended over the whole $\PMt$ by $\tdel$-equivariance, the term $r^{-\delta}$ making this function $(-\delta)$-equivariant.

Bijectivity and naturality of $\ete{u}$ are straightforward to perform.
\end{dem}

On the other hand, every function in $\finv{\PMt}{\Vt}{\rot}^{-\delta}$ can be projected over $\M$ on a function in $\finv{\PM}{\V}{\rodel}$.
\medskip
 
\begin{lem}\label{utd}
For every representation $(\V,\rodel)$ of $\GLm$ there exists a natural map
\[
\ut{d}:\finv{\PMt}{\Vt}{\rot}^{-\delta}\to\finv{\PM}{\V}{\ro_\delta}.
\]
\end{lem}
\begin{dem} Again, it is sufficient to prove the result when $\V$ is irreducible. The map $\ut{d}$ is given by 
\[
\ut{d}(g)(\pi_*w)=r^\delta (\pi_0\circ g)(w,\Eul)_{\cl{\pi_*w,r}}\pourtt g\in\finv{\PMt}{\Vt}{\rot}^{-\delta},
\]
where $\pi_0$ denotes the projection from $\Vt$ to $\V_0$. The term $r^{\delta}$ and the $\tdel$-equivariance of its arguments ensure that $\ut{d}$ is well defined, and its naturality is easy to check. 
\end{dem}

\begin{thm}\label{Des}There exists a natural map 
\[
\Des: \mathcal{D}(\ete{\E}_1,\tilde{\E}_2)\to\mathcal{D}(\E_1,E_2).
\]
\end{thm}
\begin{dem} 
Let $\mathcal{D}$ be a differential operator acting between $(-\lambda)$-equivariants sections of $\ete{\E}_1$ and $(-\mu)$-equivariant sections of $\tilde{\E}_2$. Since the two vertical arrows in the following diagram are natural

\begin{equation*}
\xymatrix@=.2em{\finv{\PMt}{\Vete_1}{\rot_1}^{-\lambda}\ar[rrrrrr]^{\mathcal{D}}
&&&&&&
\finv{\PMt}{\Vt_2}{\rot_2}^{-\mu}\ar[dddd]^{\ut{d}}\\
&&&&&&\\
&&&&&&\\
&&&&&&\\
\finv{\PM}{\V_1}{\ro_{1\lambda}}\ar[uuuu]^{\ete{u}}&&&&&&
\finv{\PM}{\V_2}{\ro_{2\mu}}}
\end{equation*}
the map $\Des$ defined by
\[
\Des(\mathcal{D})=\ut{d}\circ\mathcal{D}\circ\ete{u}
\]
is also natural. 
\end{dem}

\vspace{1em}
\begin{thm}
Let $\delta=\mu-\lambda$ be non resonant for $\V_1^*\otimes\V_2\otimes\sym^k\R^m$ given any $k\in\N$. Then there exists a natural and projectively equivariant quantization
\[
\Quant:\Co{\M}\times\mathcal{S}(E_1,
E_2)\to\mathcal{D}(\fib{E}_1,\fib{E}_2).
\]
\end{thm}
\begin{dem} 

We claim that the map $\Quantn$ defined by
\[
\Quantn =\Des\circ\taut\circ\Reln\pourtt\nabla\in\Co{\M},
\]
has the required properties.

If the complex numbers characterizing $\V_1$ end $\V_2$ are respectively given by $\lambda$ and $\mu$, then for any $k\in\N$ the complex number characterizing $\V_1^*\otimes\V_2\otimes\sym^k\R^m$ is $\delta=\mu-\lambda$. Since it is not resonant for this representation by hypothesis, $\Rel$ is well defined and in turn $\Quant$ is well defined. It is also obviously natural and projectively equivariant since it is the composition of natural and projectively equivariant maps.

The fact that $\Quantn$ preserves the principal symbol is a consequence of the fact that  
\[
\sigma(\Des(\taut g))=\sigma(\tau_\nabla(g_0))\pourtt g\in\finv{\PMt}{\Vt}{\tdel},
\]
and
\[
(\Reln f)_0=f\pourtt f\in\finv{\PM}{\V}{\rodel}.
\]
Both are direct consequences of the way we constructed $\Des$ and $\Rel$. We thus get for every $f\in\finv{\PM}{\V}{\rodel}$

\vspace{-0.5em}
\begin{eqnarray*}
\sigma (\Quantn f)&=&\sigma (\Des(\taut(\Reln f)))\\
&=&\sigma(\tau_\nabla(\Reln f)_0)\\
&=&\sigma(\tau_\nabla f)\\
&=&f.
\end{eqnarray*}
\vspace{-0.5em}
\end{dem}
\subsection{Proof of Proposition \ref{rw}}
We begin with a lemma allowing to check $\Hmu$-equivariance of a function by means of its invariance relatively to the Lie algebra $\hmu$. 

\medskip
\begin{lem}\label{lemme3}
 Let $g\in\finv{\CM}{\Vt}{}$ be a function satisfying
\begin{equation}\label{tdelinvalg}
\Lie{\fond{h}}g=-\dtdel(h)g\pourtt h\in\hmu,
\end{equation}
and
\begin{equation}\label{moinsid}
g(u\cdot\matHmu{-\Id}{0}{0}{1})=\tdel\matHmu{-\Id}{0}{0}{1}^{-1}g(u).
\end{equation}
Then $g\in\finv{\CM}{\Vt}{\tdel}$, namely satisfies
\begin{equation}\label{tdelinvgroupe}
g(u\cdot H)=\tdel(H^{-1})g(u)\pourtt H\in\Hmu.
\end{equation}
 \end{lem}
 \begin{dem}

Since equation (\ref{tdelinvalg}) is the derivative at $t=0$ of 
 \[
g(u\cdot\exp(th))=\tdel(\exp(th)^{-1})g(u)\pourtt h\in\hmu,
 \]
it implies that (\ref{tdelinvgroupe}) is satisfied for the elements or $\Hmu$ that belong to the connex component of the identity.
It is known that $\Hmu$ is connex whenever $m$ is even, and has two connex components if $m$ is odd. The result is therefore obvious if $m$ is even.

If $m$ is odd, then every point in the connex component of $\Hmu$ which does not contain the identity can be decomposed into
$H\matHmu{-\Id}{0}{0}{1}$, where $H$ belongs to the connex component of the identity.
Consequently, conditions (\ref{tdelinvalg}) and (\ref{moinsid}) imply
\begin{eqnarray*}
g(u\cdot
H\matHmu{-\Id}{0}{0}{1})&=&\tdel\matHmu{-\Id}{0}{0}{1}^{-1}g(u\cdot H)\\
&=&\tdel\matHmu{-\Id}{0}{0}{1}^{-1}\tdel(H^{-1})g(u)\\
&=&\tdel(H\matHmu{-\Id}{0}{0}{1})^{-1}g(u). 
\end{eqnarray*}
\end{dem}

We are going to prove Proposition \ref{rw} when $\V$ is an irreducible representation of $\GLm$. Indeed any $\V$ can be decomposed into a finite direct sum of irreducible components, and the lift of symbols is defined component wise. 

\medskip

Suppose that $g$ is a $\tdel$-equivariant function. From Proposition \ref{castinv}, it is an eigenvector of $\casw$ (of eigenvalue $\kappa$) if it satisfies

\begin{eqnarray*}
\kappa g&=&\sum_{q\in Q}\al_q g_q-2\sum_{q\in Q}\sum_{i=1}^m
i(\epsilon^i)\Liew{v_i} g_q.
\end{eqnarray*}
By projecting this equation on each irreducible component $\V_q$ of $\Vt$, we get
\begin{eqnarray*}
\kappa g_0&=&\alpha_0 g_0;\\
\kappa g_q&=&\al_q g_g-2\sum_{\va{q'}=\va{q}-1}i(\epsilon^i)\Liew{v_i}
g_{q'}\pourtt q\in Q,q\neq 0.
\end{eqnarray*}
Since $\delta$ is not resonant for $\V$, we know that $\al_q\neq\al_0$ for $q\neq 0$, thus these equations are equivalent to
\begin{equation}\label{vpcas}
g_q=-\frac{2}{(\al_0-\al_q)}\sum_{\va{q'}=\va{q}-1}i(\epsilon^i)\Liew{v_i}
g_{q'}\pourtt q\neq0,
\end{equation}
$g_0$ being left arbitrary. In particular, if $g\in\finv{\CM}{\Vt}{\tdel}$ is an eigenvector of $\casw$, then it is uniquely determined by its projection $g_0$ over $\V_0$. This leads us to set
\[
\Rw: \finv{\CM}{\V}{\sdel}\to\finv{\CM}{\Vt}{}
\]
as associating to any $\sdel$-equivariant $\V$ valued function $f$ the unique $\Vt$ valued function satisfying
\begin{itemize}
\item[$\bullet$]
$(\Rw(f))_0=f$;
\item[$\bullet$]$\displaystyle{(\Rw(f))_q=-\frac{2}{(\al_0-\al_q)}\sum_{\va{q'}=\va{q}-1}i(\epsilon^i)\Liew{v_i}
(\Rw(f))_{q'}\pourtt q\neq0}$.
\end{itemize}

\medskip
 As it is defined, $\Rw$ could be valued in arbitrary functions in $\finv{\CM}{\Vt}{}$. We check thereafter that it is actually valued in $\tdel$-equivariant function, as announced in the statement of Proposition \ref{rw}. More precisely, we prove that any function from 
$g\in\finv{\CM}{\Vt}{}$ satisfying (\ref{vpcas}) and for which
$g_0$ is ${\sdel}$-equivariant is necessarily $\tdel$-equivariant.

To that aim, let us set 
\[
\casw_0(g)=\sum_{q\in Q}\al_q g_q-2\sum_i i(\epsilon^i)\Liew{v_i}g.
\]
This operator coincides exactly with $\casw$ on $\tdel$-equivariant functions, and functions satisfying (\ref{vpcas}) are eigenvectors of $\casw_0$. A quite long but straightforward computation shows that
\begin{equation}\label{nono}
(\dtdel(h)+\Lie{\fond{h}})\circ
\casw_0=\casw_0\circ(\dtdel(h)+\Lie{\fond{h}})\pourtt h\in\hmu.
\end{equation} 
Together with the fact that $g$ is an eigenvector of $\casw_0$, it implies that $(\dtdel(h)+\Lie{\fond{h}})g$ is also an eigenvector of $\casw_0$, thus determined by its projection on $\V_0$. We are going to show that the latter vanishes for any $h\in\hmu$, implying that $(\dtdel(h)+\Lie{\fond{h}})g=0$ and hence Condition (\ref{tdelinvalg}).

If $h\in\go$, then the $\sdel$-equivariance of $g_0$ implies that 
\begin{eqnarray*}
[(\dtdel(h)+\Lie{\fond{h}})g]_0&=&\dsdel(h)g_0+\Lie{\fond{h}}g_0\\
&=&0.
\end{eqnarray*}

If $h\in\gun$, then $\dtdel(h)g$ is valued in $\oplus_{q\neq 0}\V_q$ and due to the $\sdel$-equivariance of $g_0$, we get $\Lie{\fond{h}}g_0=0$. Thus,
\[
[(\dtdel(h)+\Lie{\fond{h}})g]_0=0.
\] 

Moreover, it is easily verified that
\begin{equation}\label{la}
g(*\cdot\matHmu{-\Id}{0}{0}{1})=\tdel\matHmu{-\Id}{0}{0}{1}^{-1}g(*),
\end{equation}
by projecting this identity on every $\V_q$ and using the explicit form of each $g_q$ given in (\ref{vpcas}). We can thus apply Lemma \ref{lemme3}, and conclude that $\Rw$ is valued in $\finv{\CM}{\Vt}{\tdel}$.
\medskip

It remains to show the naturality of $\mathcal{R}$. Let $\phi$ be a local diffeomorphism over $\M$. The map $\mathcal{R}$ is natural if it satisfies for all Cartan connection $\w$ and all $f$ in $\finv{\CM}{\V}{\sdel}$
\begin{equation}\label{natReln}
\phipb(\Rw f)=\mathcal{R}^{\phipb\w}(\phipb f),
\end{equation}
which is the case if $\phipb(\Rw f)$ is the (unique) eigenvector of $\cas{\phipb\w}$ whose projection on $\V_0$ is equal to $\phipb f$.
On the one hand, the naturality of the Casimir operator implies
\begin{eqnarray*}
\cas{\phipb\w}(\phipb(\Rw f))&=&\phipb\casw (\Rw f)\\
&=&\al_0\phipb(\Rw f).
\end{eqnarray*}
On the other hand, we have

\begin{eqnarray*}
(\phipb \Rw f)_0&=&\phipb(\Rw f)_0\\
&=&\phipb f,
\end{eqnarray*}
which ends the proof of Proposition \ref{rw}.

\section*{Acknowledgements} I would like to warmly thank Pierre Lecomte and Pierre Mathonet for fruitful suggestions and remarks, and for all the stimulating mathematical discussion we had, together with Fabian Radoux. I also thank Bruno Teheux, R\'{e}mi Lambert and Michel Georges for checking the draft of this paper, and the National Fund of Scientific Research of Belgium for my fellowship.

\end{document}